\documentclass[fleqn,12pt]{article}
\usepackage{amsmath,amssymb,amsmath}
\usepackage{graphics}
\input epsf

\newtheorem{theorem}{Theorem}
\newtheorem{lemma}{Lemma}
\newtheorem{proposition}{Proposition}

\newtheorem{corollary}{Corollary}

\def\square{\hbox{\vrule\vbox{\hrule\phantom{o}\hrule}\vrule}}

\begin{document}
\title{Distances Between Composition Operators}
\author{Valentin Matache}
\date{ }

\maketitle

\begin{abstract}
Composition operators $C_\varphi $ induced by a selfmap $\varphi $
of some set $S$ are operators acting on a space consisting of
functions on $S$ by composition to the right with $\varphi $, that
is $C_\varphi f=f\circ \varphi $. In this paper, we consider the
Hilbert Hardy space $H^2$ on the open unit disk and find exact
formulas for distances $\| C_\varphi -C_\psi \| $ between
composition operators. The selfmaps $\varphi $ and $\psi $
involved in those formulas are constant, inner, or analytic
selfmaps of the unit disk fixing the origin.
\end{abstract}
\medskip


\smallskip

\smallskip

\section{Introduction}

Let $H^p$ denote the  Hardy space of index $p$ on the open unit
disk $\mathbb{U}$, that is the space of all functions $f$ analytic
in $\mathbb{U}$ satisfying the condition
\begin{equation}
\| f\| _p:=\sup_{0<r<1}\left( \int_\mathbb{\partial
\mathbb{U}}|f(r\zeta )|^p\, dm(\zeta )\right)^{1/p}<\infty ,
\end{equation}
where $m$ is the normalized Lebesgue measure and $p$ is fixed
$0<p<+\infty $.

It is well known that $\| \quad \| _2$ is a Hilbert norm on $H^2$
with alternative description
\begin{equation}
\|f\| _2=\sqrt{\sum_{n=0}^\infty |c_n|^2},
\end{equation}
where $\{ c_n\} $ is the sequence of Maclaurin coefficients of
$f$. The Hilbert Hardy space $H^2$ is our space of choice in this
paper.

The space $H^\infty $ is the space of all bounded analytic
functions on $\mathbb{U}$ endowed with the supremum norm $\| \quad
\| _\infty  $. It is easy to see that $H^\infty \subseteq H^p$,
$0<p<+\infty $. Another well known fact about $H^p$--functions is
the fact that, by a classical result of P. Fatou \cite[Theorem 1.3
]{Duren}, eventually extended by F. and M. Riesz, those functions
have nontangential limits a.e. on $\partial \mathbb{U}$. The
nontangential limit function of any $f$ in $H^p$ will be denoted
by the same symbol as the function itself. It is known that it is
an $L^p_{\partial \mathbb{U}}$--function and
\begin{equation}
\|f \| _p=\left(\int_\mathbb{\partial \mathbb{U}}|f(\zeta )|^p\,
dm(\zeta )\right)^{1/p}\qquad f\in H^p.
\end{equation}
An analytic selfmap of $\mathbb{U}$ is called an inner function if
it has unimodular nontangential limits a.e. on $\partial
\mathbb{U}$.

For each analytic selfmap $\varphi $ of $\mathbb{U}$ \textit{the
composition operator of symbol} $\varphi $ is the following
operator
\begin{equation}
C_\varphi f=f\circ \varphi \qquad f\in H^2.
\end{equation}
Such operators are bounded, as a consequence of Littlewood's
Subordination Principle, \cite[Theorem 1.7]{Duren}, saying that
composition operators whose symbol fixes the origin are
contractions. If $\varphi $ is a conformal automorphism of
$\mathbb{U}$, we call $C_\varphi $  an automorphic composition
operator.

The numerical range of a Hilbert space operator $T$ is the set
$W(T)=$

\noindent$\{ <Tf,f>: \| f\| =1\} $. It is well known that
numerical ranges are convex subsets of the complex plane whose
closure contains the spectrum of the given operator, \cite[Chapter
22]{Halmos}. The quantity $w(T)=\sup \{ |<Tf,f>|: \|f\| =1\} $ is
called the numerical radius of the operator $T$.

A Hilbert space operator $T$ that satisfies an equation of the
form $ T^2+\lambda T+\mu I=0 $ where $I$ is the identity operator
and $\lambda $ and $\mu $ are constants is called a quadratic
operator.

In section 2 of this paper we are able to complete the description
of numerical ranges of quadratic, automorphic composition
operators by showing that the numerical range of such an operator
is open, unless the operator is $C_z=I$ or $C_{(-z)}$. The
description of the closure of the aforementioned numerical range
was obtained in \cite{Abdollahi} and \cite{Bourdon}.

Quadratic Hilbert space operators are known to have elliptical
numerical ranges, \cite[Theorem 2.1]{Tso}. We apply the
aforementioned result and the theorems characterizing the
numerical ranges of quadratic composition operators to calculate
distances between composition operators. The symbols of the
composition operators involved in these distance computations are
constant or inner. We use different methods to calculate the
distance $\| C_0-C_\varphi \| $ when $\varphi (0)=0$. We show that
the distance above equals some Hardy norm $\| \varphi \| _p$, of
$\varphi $, for some $2\leq p\leq +\infty $. We are able to show
that any  $2\leq p\leq +\infty $ works if and only if $\varphi $
is a scalar multiple of an inner function. Otherwise, we show the
choice of $p$ is unique, and that $p=2$ if and only if $\varphi $
is orthogonal to the set $\{ \varphi ^2, \varphi ^3, \dots \} $.
 These distance computations are in
section 3.

\section{Description of a Numerical Range}

By the spectral mapping theorem \cite[Chapter 9]{Halmos}, the
spectrum of a quadratic operator can consist of at most 2 points.
Quadratic operators with spectrum consisting of two points are
known to have elliptical numerical ranges. More exactly, the
following is proved in \cite[Theorem 2.1]{Tso}.

\begin{theorem}
The numerical range of a quadratic operator having spectrum
consisting of the two distinct points $a$ and $b$ is an open or a
closed elliptical disk, possibly degenerate, {\rm (}that is,
reduced to its focal axis{\rm )}. The major axis of the disk has
length $\| T-aI\| $ and the length of the minor axis is $\sqrt{(\|
T-aI\| ^2-|a-b|^2)}$. The elliptical disk is closed if and only if
$T$ attains its norm or equivalently, if and only if it attains
its numerical radius.
\end{theorem}

Above, the statement that $T$ attains its norm, respectively its
numerical radius, means that there is some norm--one vector $f$ so
that $\| T\| =\| Tf\|$, respectively $| <Tf,f>|=w(T)$.

It is easy to determine the quadratic automorphic composition
operators. Indeed, by \cite{Nordgren01}, the only automorphic
composition operators with spectrum consisting of at most 2 points
are the identity operator and the composition operators whose
automorphic symbol $\varphi $ fixes a point in $\mathbb{U}$ and is
conformally conjugated to $-z$, that is $\varphi $ should be of
the form
\[
\varphi (z)=\alpha _p(z)=\frac{p-z}{1-\overline{p}z}\qquad z\in
\mathbb{U},
\]
where $p$ is any fixed point in $\mathbb{U}$. Visibly such
operators are quadratic because $\alpha _p\circ \alpha _p(z)=z$,
$z\in \mathbb{U}$ and hence $C^2_{\alpha _p}=I$.

The closure of $W(C_{\alpha _p})$ was characterized by the authors
of \cite{Bourdon}. They showed it is a closed elliptical disc of
foci $-1$ and $1$. That disk is reduced to its focal axis if and
only if $p=0$. The authors of \cite{Bourdon} gave a formula for
the length of the major axis of that disk. That formula is hard to
use in practical problems. Therefore, very recently, the author of
\cite{Abdollahi} found the following practical formula for the
length of the aforementioned major axis.
\begin{theorem}
For each $p\in \mathbb{U}$, the length of the major axis of the
closure of $W(C_{\alpha _p})$  is $2/\sqrt{1-|p|^2}$.
\end{theorem}

Except the case when $p=0$, it is not known if $W(C_{\alpha _p})$
is open or closed. We prove in the sequel that it is open.
According to Theorem 1, we should check if $C_{\alpha _p}$ attains
its norm or not.

Norms of composition operators are not easy to calculate. The
commonest things known about them are the following.
\begin{equation}
\frac{1}{\sqrt{1-|\varphi (0)|^2}}\leq \| C_\varphi \| \leq \sqrt{
\frac{1+|\varphi (0)|}{1-|\varphi (0)|}}
\end{equation}
Visibly the lower and upper bound coincide to $1$ if $\varphi $
fixes the origin. If $\varphi (0)\neq 0 $, then the lower bound is
attained if and only if the symbol is constant
 \cite[Theorem 4]{Pokorny-Shapiro}, whereas the upper bound
is attained if and only if the symbol is an inner function
\cite[Theorem 5.2]{Shapiro-Monatschefte}. We are now ready to
prove that the numerical range of a quadratic, automorphic
composition operator, other than $C_{(-z)}$ or $C_z$ is an open
elliptical disk. We obtain this result, as a consequence of the
following.

\begin{proposition}
A composition operator having inner symbol $\varphi $ attains its
norm if and only if $\varphi (0)=0$.
\end{proposition}
\medskip

\noindent \textit{Proof.} Let
\[
P(z,u)=\Re \frac{u+z}{u-z}\qquad u\in \partial \mathbb{U},z\in
\mathbb{U}
\]
be the usual Poisson kernel.

Let $\varphi $ be an inner function. The following formula is
established in \cite{Nordgren01}
\begin{equation}
\int_{\partial \mathbb{U}}|f\circ \varphi (u)|^2\,
dm(u)=\int_{\partial \mathbb{U}}|f(u)|^2P(\varphi (0),u)\,
dm(u)\qquad f\in H^2.
\end{equation}

An immediate consequence is the fact that composition operators
whose symbols are inner functions fixing the origin are
isometries.

It is well known and easy to prove that
\begin{equation}
\frac{1+|z|}{1-|z|}\geq P(z,u) \qquad u\in \partial \mathbb{U},
z\in \mathbb{U}.
\end{equation}
Given (6), if $f\in H^2 $ has norm 1, the relation $\| C_\varphi
f\| _2=\| C_\varphi \|$ is equivalent to
\[
\int_{\partial \mathbb{U}}|f(u)|^2\left(\frac{1+|\varphi
(0)|}{1-|\varphi (0) |}-P(\varphi (0),u)\right)\, dm(u)=0.
\]
By (7), it follows that
\[
|f(u)|^2\left(\frac{1+|\varphi (0)|}{1-|\varphi (0) |}-P(\varphi
(0),u)\right)=0\qquad {\rm a.e.}
\]
Since $\| f\| _2=1$, we deduce
\[
\frac{1+|\varphi (0)|}{1-|\varphi (0) |}=P(\varphi (0),u)\qquad
{\rm a.e.},
\]
a condition that is satisfied if and only if $\varphi (0)=0$.
Since composition operators with inner symbol fixing the origin
are isometric, they obviously attain their norm. \hfill \square
\bigskip

\begin{corollary}
If $p\neq 0$, $W(C_{\alpha _p})$ is the open elliptical disk of
foci $\pm 1$ and major axis of length $2/\sqrt{1-|p|^2}$.
\end{corollary}
\medskip

\noindent \textit{Proof.} This is a direct consequence of Theorem
1, Proposition 1, and the computation of the length of the major
axis contained by Theorem 2.\hfill \square \bigskip

\section{Distance Computations}

Composition operators of constant symbol $C_p$, $p\in \mathbb{U}$
are obviously idempotent and hence quadratic. As one can see, we
denote by $p$ both the function on $\mathbb{U}$ constantly equal
to $p$ and the complex number $p$ itself. Since composition
operators having constant symbol are rank-one operators, they are
compact on the infinite dimensional space $H^2$. Thus their
spectrum contains $0$. The evident relation $C_\varphi 1=1$, valid
for any composition operator, shows that the spectrum of $C_p$
consists of $0$ and $1$. Then the following theorem originally
proved in \cite{Matache_Nranges}, (see also
\cite{Bourdon_Shapiro_0}), becomes very easy to obtain as a
consequence of Theorem 1, the formula for the norm of composition
operators of constant symbol, and the fact that composition
operators of constant symbol attain their norms, (which is not
hard to prove).
\begin{theorem}
The numerical range of $C_p$ is  the closed elliptical disk of
foci $0$ and $1$ and major axis $\frac{1}{\sqrt{1-|p|^2}}$. The
disk is reduced to its focal axis if and only if $p=0$.
\end{theorem}

From now on, our basic problem will be finding  the norm of a
difference of two distinct composition operators.

First,
 recall that $H^2$ is a reproducing kernel Hilbert space,
(see \cite[Chapter 4]{Halmos} for the basics on this kind of
spaces), that is the functions $k_p(z)=1/(1-\overline{p}z)$,
(called the kernel--functions of $H^2$) have the property below,
known as "the
 reproducing property"
\[
f(p)=<f,k_p>\qquad p\in \mathbb{U}, \, f\in H^2.
\]

By the "reproducing property", for any $p_1,p_2\in \mathbb{U}$,
the following holds
\begin{equation}
 \| k_{p_1}-k_{p_2}\|
=\sqrt{\frac{1}{{1-|p_1|^2}}+\frac{1}{{1-|p_2|^2}}-2\Re
\frac{1}{{1-\overline{p_1}p_2}}}.
\end{equation}

Next, let us introduce more terminology. For an analytic map $\psi
$ on $\mathbb{U}$, $M_\psi $ denotes the multiplication operator
of symbol $\psi $, that is the operator
\[
M_\psi f=\psi f\qquad f\in H^2.
\]
It is well known that $M_\psi $ is bounded on $H^2$ if and only if
$\psi \in H^\infty $ and, in that case, $\| M_\psi \| =\| \psi \|
_\infty  $. For an analytic map $\psi $ on $\mathbb{U}$ and an
analytic selfmap $\varphi $ of $\mathbb{U}$, the operator $
T_{\psi ,\varphi }=M_\psi C_\varphi  $ is called the weighted
composition operator of symbols $\psi $ and $\varphi $. Clearly
such operators are bounded if the first symbol is bounded.

 Denote $H^2_0=zH^2=H^2\ominus
\mathbb{C}$. A theorem that will be extensively cited in the
sequel is
 the following, \cite[Theorem 5.1]{Shapiro-Monatschefte}.

\begin{theorem}
Let $\varphi $ be an analytic selfmap of $\mathbb{U}$ that fixes
the origin. Then $\| C_\varphi |H^2_0\| =1$ if and only if
$\varphi $ is inner.
\end{theorem}

We give a new proof of Theorem 4, based on the following result
which appears in \cite{Shapiro-Monatschefte} with a proof
independent of that of Theorem 4 and was recently given an
alternative short proof in \cite{Matache-Inner}.\bigskip

\textit{The function $\varphi :\mathbb{U}\rightarrow \mathbb{U}$
is inner if and only if
\[
 \| C_\varphi \| _e=\sqrt{\frac{1+|\varphi (0)|}{1-|\varphi (0)|}},
\]
where $\| C_\varphi \| _e$ denotes the essential norm of}
$C_\varphi $.
\bigskip

Before giving the announced new proof to Theorem 4 we record in a
lemma a fact that is probably known as folklore and equally easy
to prove.

\begin{lemma}
If $T$ is a Hilbert--space operator on the space $H$ and $L$ is a
closed invariant subspace of $T$, then $\| T|L\| _e\leq \| T\|
_e$.
\end{lemma}
\medskip

\noindent \textit{Proof.} Indeed, if $P$ is the orthogonal
projection of $H$ onto $L$ and $K$ any compact operator on $H$,
one can write
\[
\| T+K\| \geq \| P(T+K)|L\| =\| (T|L)+(PK|L)\| \geq \| T|L\| _e
\]
since $PK|L$ is a compact operator on $L$.\hfill \square \bigskip

\noindent \textit{Proof of Theorem $4$.} As we already noted, if
$\varphi $ is an inner function fixing the origin, then $C_\varphi
$ is isometric, hence $\| C_\varphi |H^2_0\| =1$. The delicate
part is the converse implication. Note that, if $\varphi (0)=0$,
then $C^*_\varphi 1=1$, thus $H^2_0$ is an invariant subspace of
$C_\varphi $. Therefore, if $\| C_\varphi |H^2_0\| =1$, then we
distinguish between two cases. If $\| C_\varphi |H^2_0\| _e=1$,
then by Lemma 1, $\| C_\varphi \| _e=1$, since $C_\varphi $ is a
contraction. Thus,  $\| C_\varphi \| _e=\sqrt{(1+|\varphi
(0)|)/(1-|\varphi (0)|)}$ and hence $\varphi $ must be inner. The
other case is when $\| C_\varphi |H^2_0\| _e<1=\| C_\varphi
|H^2_0\| $. In this second case, it is very easy to prove that
$C_\varphi |H^2_0$ is a norm--attaining operator,
\cite[Proposition 2.2]{Hammond}. Hence, one can consider a
norm--one $f\in H^2$ so that
\[
1=\| C_\varphi (zf(z))\| ^2_2=\| \varphi f\circ \varphi \|
^2_2\leq \| f\circ \varphi \| ^2_2\leq \| f\| ^2_2=1.
\]
One gets that
\[
\int_{\partial \mathbb{U}}|f\circ \varphi (\zeta )|^2(1-|\varphi
(\zeta ) |^2)\, dm(\zeta )=0.
\]
Since, $|\varphi (\zeta )|\leq 1$ a.e., one deduces $|f\circ
\varphi (\zeta )|^2(1-|\varphi (\zeta ) |^2)=0$ a.e., which
implies that $|\varphi (\zeta )|= 1$ a.e., that is $\varphi $ is
inner.\hfill \square
\bigskip

 In the next proposition we prove that the norm of the
restriction to $H^2_0$ of a composition operator is always equal
to the distance between itself and certain composition operators
having constant symbols. We also show that the distance between a
composition operator $C_\varphi $ and the orthogonal projection
$C_0$ onto the subspace $\mathbb{C}$ of constant functions equals
the norm of the weighted composition operator of identical
symbols, $T_{\varphi, \varphi }$.

\begin{proposition}
For each composition operator on $H^2$ the following relations
hold
\begin{equation}
\| \varphi \|_2\leq \| C_\varphi -C_0\| =\| C_\varphi |H^2_0\|=\|
T_{\varphi ,\varphi }\|  \leq \| \varphi \|_\infty \| C_\varphi
\|.
\end{equation}
For any  analytic selfmap $\varphi $  of $\mathbb{U}$
\begin{equation}
\| C_\varphi -C_{\varphi (0)}\| \leq \| T_{\varphi ,\varphi } \|
\quad \mbox{hence}\quad \| C_\varphi -C_{\varphi (0)}\| \leq \|
C_\varphi \| .
\end{equation}
The second inequality in {\rm (10)} is strict if $\| \varphi \|
_\infty <1$ or if $\varphi $ is a non--inner function fixing the
origin. If $\varphi $ is inner, then
\begin{equation}
\| C_\varphi -C_{\varphi (0)}\| =\| C_\varphi -C_0\| =\| C_\varphi
|H^2_0\|=\| T_{\varphi ,\varphi }\|  = \| C_\varphi \| .
\end{equation}

\end{proposition}

\medskip

\noindent \textit{Proof.} Clearly $\| C_\varphi -C_0 \| =\|
C_{\varphi }|H^2_0\| $, since  $(C_\varphi -C_0)1=0$ and
$C_0|H^2_0=0$. Each function $g$ in $H^2_0$ factors as
$g(z)=zf(z)$ where $f$ is an $H^2$--function of norm equal to the
norm of $g$, and obviously $C_\varphi g=T_{\varphi ,\varphi }f$.
Hence $\| C_\varphi -C_0\| =\| C_\varphi |H^2_0\|=\| T_{\varphi
,\varphi }\| $. Obviously
\[
\| T_{\varphi ,\varphi }\| \leq \| M_\varphi \| \| C_\varphi \|
=\| \varphi \| _\infty \| C_\varphi \|,
\]
which proves the upper estimate in (9). To prove the lower
estimate, take the norm--one function $z$ and note that
$(C_\varphi |H^2_0)z=\varphi $.

Note now that the range $R(C_\varphi -C_{\varphi (0)})$ of
$C_\varphi -C_{\varphi (0)}$ is contained in $H^2_0$. Therefore,
 $C_\varphi -C_{\varphi (0)}$ leaves $H^2_0$ invariant. Any
difference of composition operators transforms $1$ into $0$, thus
$H^2_0$, actually reduces $C_\varphi -C_{\varphi (0)}$, and
$C_\varphi -C_{\varphi (0)}=0\oplus ((C_\varphi -C_{\varphi (0)
})|H^2_0)$. Keeping this in mind, observe that
\[
\| (C_\varphi -C_{\varphi (0)})f\| _2^2=\sum_{n=1}^\infty| <f\circ
\varphi -f(\varphi (0)),z^n>| ^2=
\]
\[
\sum_{n=1}^\infty| <f\circ \varphi ,z^n>| ^2=\| C_\varphi f\|
_2^2-|f(\varphi (0)|^2\leq \| C_\varphi f\|^2_2\qquad f\in H^2.
\]

Hence $\| C_\varphi -C_{\varphi (0)}\| \leq \| C_\varphi \| $.
Also, substitute $f$ by $zf(z)$ above getting
\[
\| (C_\varphi -C_{\varphi (0)})(zf(z))\| _2^2=\| C_\varphi
(zf(z))\| _2^2-|\varphi (0)f(\varphi (0))|^2\leq \| T_{\varphi
,\varphi }f\| _2^2\qquad f\in H^2.
\]
Since $C_\varphi -C_{\varphi (0)}=0\oplus ((C_\varphi -C_{\varphi
(0) })|H^2_0)$, one gets $ \| C_\varphi -C_{\varphi (0)}\| \leq \|
T_{\varphi ,\varphi } \| $ which proves (10) and, by (9), also
shows that the second inequality in (10) is a strict inequality if
$\| \varphi \| _\infty <1$. The fact that the same inequality is
strict if $\varphi (0)=0$ and $\varphi $ is not an inner function
is a direct consequence of Theorem 4 and (9).

As we noted,  $\| C_\varphi \| _e=\sqrt{(1+|\varphi
(0)|)/(1-|\varphi (0)|)}=\| C_\varphi \| $ if $\varphi $ is inner.
On the other hand $\| C_\varphi \| _e\leq \| C_\varphi -C_{\varphi
(0)}\| \leq \| C_\varphi \| $, hence (11) holds if $\varphi $ is
inner.\hfill \square
\bigskip

\begin{corollary}
For any analytic selfmap $\varphi $ of $\mathbb{U}$  the equality
\begin{equation}
 \| C_\varphi -C_0\| =\sqrt{\frac{1+|\varphi
(0)|}{1-|\varphi (0)|}}
\end{equation}
holds if and only if $\varphi $ is inner. Hence composition
operators having inner symbol attain their essential norm.
\end{corollary}
\medskip

\noindent \textit{Proof.} If $\varphi $ is inner, then (12) is a
consequence of (11) and the formula for the norm of a composition
operator of inner symbol.  Conversely, if the equality in (12)
holds, then by (9), one has $\sqrt{(1+|\varphi (0)|)/(1-|\varphi
(0)|)}=\| C_\varphi |H^2_0\| \leq \| C_\varphi \| \leq
\sqrt{(1+|\varphi (0)|)/(1-|\varphi (0)|)}$, and hence $\varphi $
is an inner function if $\varphi (0)\neq 0$. If $\varphi (0)=0$
and (12) holds, then $\varphi $ is inner, by (9) and Theorem 4.
\hfill \square \bigskip

The following proposition contains two very simple
distance--formulas.

\begin{proposition}
Let $p_1,p_2\in \mathbb{U}$, $\lambda ,\mu \in
\overline{\mathbb{U}}$, and $\varphi $ be an inner function,
fixing the origin. Then
\begin{equation}
\| C_{p_1} -C_{p_2}\| =
\sqrt{\frac{1}{{1-|p_1|^2}}+\frac{1}{{1-|p_2|^2}}-2\Re
\frac{1}{{1-\overline{p_1}p_2}}}
\end{equation}
and
\begin{equation}
\| C_{\lambda \varphi }-C_{\mu \varphi }\| =\sup \{ |\lambda
^n-\mu ^n| : n=1,2,3,\dots \} .
\end{equation}
\end{proposition}
\medskip

\noindent \textit{Proof.} Since, for each $f\in H^2$, $\| (C_{p_1}
-C_{p_2})f\| _2=|f(p_1)-f(p_2)|=$

\noindent $|<k_{p_1}-k_{p_2},f>|$, (13) is a consequence of (8).

To prove (14), note that $\| C_{\lambda \varphi }-C_{\mu \varphi
}\|=\| C_\varphi (C_{\lambda z}-C_{\mu z})\| =\| C_{\lambda
z}-C_{\mu z}\|$. The operator $(C_{\lambda z}-C_{\mu z})$ has a
diagonal matrix in the standard Hilbert base $B=\{
1,z,z^2,z^3,\dots, z^n, \dots \} $ of $H^2$. The diagonal entries
are $\{ 0, \lambda -\mu , \lambda ^2-\mu ^2, \lambda ^3-\mu
^3,\dots \} $ and hence, (14) follows, by the well known formula
for the norm of a diagonal operator, \cite[Chapter 7]{Halmos}.
\hfill \square \bigskip

\begin{corollary}
If $\varphi $ is an inner function fixing the origin and $\lambda
$, $\mu \in \partial \mathbb{U}$, $\lambda \neq \mu $, then $\|
C_{\lambda \varphi }-C_{\mu \varphi }\| =2$ if $\lambda /\mu $ is
a  root of unity of even order or if $\lambda /\mu $ is not a root
of unity. If $\lambda /\mu $ is a  root of unity of odd order $k$,
then $\| C_{\lambda \varphi }-C_{\mu \varphi }\| =|1-e ^{\pi
(k-1)/k}|$.
\end{corollary}
\medskip

\noindent \textit{Proof.} If $\lambda /\mu $ is not a root of
unity, then the set of powers of $\lambda /\mu $ is a dense subset
of $\partial \mathbb{U}$. This, the identity $|\lambda ^n-\mu
^n|=|\lambda ^n/\mu ^n-1|$, and formula (14) imply that $\|
C_{\lambda \varphi }-C_{\mu \varphi }\| =2$, when $\lambda /\mu $
is not a root of $1$. If $\lambda /\mu $ is a root of unity, the
formulas for $\| C_{\lambda \varphi }-C_{\mu \varphi }\| $ in this
corollary are direct consequences of formula (14) and the
geometric representation of the roots of unity as vertices of a
regular polygon inscribed in $\partial \mathbb{U}$, having a
vertex at $1$.\hfill \square \bigskip

\begin{corollary}
Denote $\varphi ^{[n]}=\varphi \circ \dots \varphi $, $n$ times.
If $\varphi $ has a fixed point in $\mathbb{U}$ and is not an
inner function, then the second inequality in {\rm (10)} is strict
for all $\varphi ^{[n]}$ starting some $n$.
\end{corollary}
\medskip

\noindent \textit{Proof.} Choose any $\varphi $, non--inner, with
a fixed point $p\in \mathbb{U}$. If $p=0$, then  one has $\|
C_{\varphi ^{[n]}}-C_{\varphi ^{[n]}(0)}\| < \| C_{\varphi
^{[n]}}\|$  for all $n$, by (9) and Theorem 4.

If $p\neq 0$, note that, by Schwarz's lemma in classical complex
analysis,  $\varphi ^{[n]}\to p$ uniformly on compacts, (see also
\cite{Shapiro1}). This combined with formula (13) proves that $ \|
C_{\varphi ^{[n]}(0)}- C_p\| \to 0$. On the other hand, it is
shown in \cite{Matache_Convergence} that $\| C_{\varphi
^{[n]}}-C_p\| \to 0$. Thus, if, arguing by contradiction, one
assumes $ \| C_{\varphi ^{[n]}}-C_{\varphi ^{[n]}(0)}\| \geq \|
C_{\varphi ^{[n]}}\|$ for infinitely many values of $n$, one gets
the contradiction $0\geq \| C_p\|  $. \hfill \square
\bigskip

The next group of distance--formulas we prove are immediate
consequences of Theorems 1, 2, 3, and of the fact that an inner
function fixing the origin induces an isometric composition
operator.

\begin{proposition}
If $\varphi $ is an inner function fixing $0$ and $p$ a constant
in $\mathbb{U}$ then
\begin{equation}
\| C_\varphi -C_p\| =\frac{1}{\sqrt{1-|p|^2}}
\end{equation}
and
\begin{equation}
\| C_{\alpha _p\circ \varphi }\pm C_{\varphi }\|
=\frac{2}{\sqrt{1-|p|^2}}.
\end{equation}
If $\varphi $ is any inner function, then
\begin{equation}
\| C_{\alpha _{\varphi (0)}\circ \varphi } \pm C_\varphi
\|=\frac{2}{\sqrt{1-|\varphi (0)|^2}}.
\end{equation}

\end{proposition}
\medskip

\noindent \textit{Proof.} Let $\varphi $ be  an inner function
fixing the origin  and note that
\[
\| C_p-I\| =\| C_\varphi (C_p-I))\| =\| C_\varphi-C_p\| .
\]
By Theorem 1, $\| C_p-I\| $ equals the major axis of the ellipse
in Theorem 3, hence (15) holds.

Given that $C_\varphi $ is isometric, if $\varphi $ is inner and
$\varphi (0)=0$,
\[
\| C_{\alpha _p\circ \varphi }\pm C_\varphi \| =\| C_\varphi
(C_{\alpha _p}\pm I)\| =\| C_{\alpha _p}\pm I \|.
\]
By Theorems 1 and 2, $\| C_{\alpha _p}\pm I \|=2/\sqrt{1-|p|^2}$,
hence (16) holds.  For an arbitrary inner function $\varphi $,
consider the inner function fixing the origin $\alpha _{\varphi
(0)}\circ \varphi $ and let $p=\varphi (0)$.  Applying (16), one
gets (17). \hfill \square
\bigskip

Consider an analytic selfmap $\varphi $ of $\mathbb{U}$ with the
property $\varphi (0)=0$. It is easy to see that the space
$\mathbb{C}$ of constant functions reduces $C_\varphi $. Indeed,
on one hand, the evident relation $C_\varphi 1=1$ shows that
$\mathbb{C}$ is left invariant by any composition operator, on the
other, as was noted in the proof of Theorem 4, $C^*_\varphi 1=1$,
so $C^*_\varphi $ leaves $\mathbb{C}$ invariant. We study in the
following the quantity $\| C_\varphi -C_0\| =\| C_\varphi |H^2_0\|
=\| T_{\varphi ,\varphi }\| $.

\begin{theorem}
If $\varphi $ is not constant and $\varphi (0)=0$, then the norm
of the restriction $C_{\varphi }|H^2_0 $ satisfies the estimates
\begin{equation}
\| \varphi \| _2\leq \| C_\varphi |H^2_0\| \leq \| \varphi \|
_\infty
\end{equation}
and the following are equivalent.\smallskip

{\rm (i)} $\frac{1}{\| \varphi \| _\infty }\varphi $ is an inner
function.\smallskip

{\rm (ii)} $\| \varphi \| _\infty =\| C_\varphi |H^2_0\|
$\smallskip

{\rm (iii)} $\| \varphi \| _2=\| \varphi \| _\infty =\| C_\varphi
|H^2_0\|  $
\end{theorem}

\medskip

\noindent \textit{Proof.} The estimates in (18) are direct
consequences of (9) and the fact that $\| C_\varphi \| $ is a
contraction when $\varphi (0)=0$.

Let us note that the situation $\| \varphi \| _2=\| \varphi \|
_\infty =\| C_\varphi |H^2_0\| $ occurs if and only if $\varphi $
is a scalar multiple of an inner function, that is if and only if
there exist $\lambda $ in  the closure of ${\mathbb{U}}$ and an
inner function $\phi $ such that $\varphi =\lambda \phi $. The
fact that, if $\varphi $ has the form above, then $\| \varphi \|
_2=\| \varphi \| _\infty =\| C_\varphi |H^2_0\| =|\lambda |$  is
immediate by (18). The converse is a consequence of the fact that
obviously $ |\varphi (e^{i\theta })|\leq \| \varphi \| _\infty $
a.e. and if $\| \varphi \| _2=\| \varphi \| _\infty $, then $
\frac{1}{2\pi }\int_{-\pi}^\pi \left(\| \varphi \| ^2_\infty -|
\varphi (e^{i\theta })| ^2\right)\, d\theta =0 $ so $|\varphi
(e^{i\theta })|=\| \varphi \| _\infty $ a.e., that is $\phi
:=\varphi /\| \varphi \| _\infty $ is inner, so setting $\lambda
:=\| \varphi \| _\infty $, one has the desired representation
$\varphi =\lambda \phi $. Visibly, $\varphi $ is a scalar multiple
of an inner function if and only if (i) holds. Thus
(i)$\Longleftrightarrow $(iii). Clearly (iii)$\Longrightarrow
$(ii). To finish the proof we show that  (ii) fails if (i) fails.

The fact that, if $\varphi $ is not a scalar multiple of an inner
function then, the upper estimate in (18) fails to be an equality
 is an immediate consequence
of Theorem 4.

Indeed, let us consider the connection between the norms $\|
C_{\phi }|H^2_0 \|$ and $\| C_\varphi |H^2_0\|$, where $\varphi
=\lambda \phi $, $\lambda =\| \varphi \| _\infty $. If one denotes
by $D$ the diagonal operator on $H^2_0$ having diagonal $\{
\lambda , \lambda ^2, \dots , \lambda ^n, \dots \} $ with respect
to the standard basis $\{ z , z ^2, \dots , z ^n, \dots \} $ of
$H^2_0$ then one has
$
C_\varphi |H^2_0=(C_{\phi }|H^2_0 )D
$
and since, $\| D\| =|\lambda |$, one gets
$
 \| C_\varphi
|H^2_0\|\leq \| \varphi \| _\infty \| C_{\phi }|H^2_0 \| ,
$
which leads to
$
 \| C_\varphi
|H^2_0\|<\| \varphi \| _\infty
$
if $\phi $ is not inner.\hfill \square \bigskip

Note that, if $\varphi (0)=0$ and $\varphi $ is a scalar multiple
of an inner function, then
\[
\| C_\varphi |H^2_0\| =\| \varphi \| _p\qquad 2\leq p\leq \infty .
\]
Relative to that, we prove in the following that $\| C_0-C_\varphi
\|  $ is always equal to some $H^p$--norm of $\varphi $, the
choice of $p$, being unique if $\varphi $ is not a scalar multiple
of an inner function.
\begin{theorem}\bigskip
For each $\varphi $ there is a $p(\varphi )$, $2\leq p(\varphi
)\leq \infty $ such that
\begin{equation}
\| C_\varphi |H^2_0\| =\| \varphi \| _{p(\varphi )}
\end{equation}
and $p(\varphi )$ is finite and uniquely determined unless
$\frac{1}{\| \varphi \| _\infty }\varphi $ is an inner function.
\end{theorem}

\medskip

\noindent \textit{Proof.} The existence of a number $p(\varphi )$
satisfying (19) is a direct consequence of estimate (18), the
continuity of the map
\[
p\rightarrow \| \varphi \| _p\qquad 2\leq p<\infty,
\]
and the fact that, $\lim_{p\rightarrow \infty }\| \varphi \| _p=\|
\varphi \| _\infty  $, ~\cite[pp 70]{Rudin}. As we noted above,
any value $p(\varphi )$, $2\leq p(\varphi )\leq \infty $ satisfies
(19) if $\frac{1}{\| \varphi \| _\infty }\varphi $ is an inner
function.

Assume now that $\frac{1}{\| \varphi \| _\infty }\varphi $ is not
an inner function, then $\| C_\varphi |H^2_0\| <\| \varphi \|
_\infty $, by Theorem 4, so a number $p(\varphi )$ satisfying (19)
must be finite. To show that there is only one such number,
assume, arguing by contradiction, that there exist $2\leq
p<r<\infty $ such that
\[
\| \varphi \| _p=\| \varphi \| _r=\| C_\varphi |H^2_0\| .
\]
Applying H\"older's inequality to the functions $| \varphi | ^p $
and $1$, one can write
\[
\int_{\partial \mathbb{U}}|\varphi (\zeta )|^p\, dm(\zeta )\leq
\left( \int_{\partial \mathbb{U}}|\varphi (\zeta )|^r\, dm(\zeta
)\right) ^{\frac{p}{r}}\qquad \mbox{i.e.}\qquad \| \varphi \|
_p\leq \| \varphi \| _r.
\]

But H\"older's inequality above is an equality, under our
assumptions, so there must be a constant $c>0$ so that $| \varphi
| ^p=c$ a.e., (\cite{Rudin}, comments following Theorem 3.5), that
is $\frac{1}{\| \varphi \| _\infty }\varphi $ is an inner
function, a contradiction.

\hfill \square

\medskip

In the case when $\frac{1}{\| \varphi \| _\infty }\varphi $ is not
inner we describe in the sequel when $p(\varphi )=2$.\medskip

\begin{theorem}
Let $\varphi $ be an analytic selfmap of $\mathbb{U}$ fixing the
origin. Then the following are equivalent.
\begin{equation}
\| C_\varphi |H^2_0\| =\| \varphi \| _2
\end{equation}
\begin{equation}
C^*_\varphi C_\varphi (z)=C^*_\varphi (\varphi )=\| \varphi \|
^2_2z
\end{equation}
\begin{equation}
\mbox{The coordinate function $z$ is an eigenfunction of
}C^*_{\varphi }C_\varphi \, .
\end{equation}
\begin{equation}
<\varphi ,\varphi ^n>=0\qquad n\geq 2
\end{equation}

\end{theorem}
\medskip

\noindent \textit{Proof.} Assume $ \| C_\varphi |H^2_0\| =\|
\varphi \| _2$. Applying the Cauchy-Schwartz inequality one gets
\[
\| C_\varphi |H^2_0\| ^2=\| \varphi \| ^2_2=|<C^*_\varphi
C_\varphi (z),z>|\leq \| C^*_\varphi C_\varphi (z)\| \| z\| \leq
\| \varphi \| ^2_2.
\]
The Cauchy-Schwartz inequality being an equality if and only if
the vectors involved in it are colinear, we get that $ C^*_\varphi
C_\varphi (z)=\lambda z $ for some scalar $\lambda $, which is
easy to determine, since
\[
\lambda =<\lambda z,z>=<C^*_\varphi C_\varphi (z),z>=\| \varphi \|
^2_2.
\]
We established (20)$\Rightarrow $(21). Obviously (21)$\Rightarrow
$(22).

If (22) holds, that is, if $C^*_\varphi C_\varphi (z)=\lambda z$
for some scalar $\lambda $, then
\[
<\varphi ,\varphi ^n>=<C^*_\varphi (\varphi ),z^n>=<\lambda
z,z^n>=0, \qquad n=0,2,3,4, \dots
\]
Hence (22)$\Rightarrow $(23).

To finish, we show (23)$\Rightarrow $(20).  Assume that
\[
\varphi \perp \varphi ^n\qquad n=0,2,3,4,\dots
\]
For any polynomial $p$ in $H^2_0$, we have
\begin{equation}
\| C_\varphi p\| ^2_2 \leq \| \varphi \| _2^2\|p\|_2^2.
\end{equation}
Inequality (24) is evident if the degree of $p$ is $1$. Arguing by
induction, assume it is true for all polynomials of degree at most
$n-1$. Consider an arbitrary polynomial $p(z) =
\sum_{j=1}^nc_jz^j$ of degree $n$ and note that
\[
p\circ \varphi (z)=c_1\varphi (z)+\varphi (z)q\circ \varphi
(z)\quad \mbox{where}\quad q(z)=\sum_{j=2}^n c_jz^{j-1}.
\]
Since $\varphi \perp \varphi (q\circ \varphi )$ and the degree of
$q$ is $n-1$, one gets
\[
\| C_\varphi p\| ^2_2 =|c_1|^2\| \varphi \| ^2_2+\| \varphi
(q\circ \varphi )\| _2^2\leq
\]
\[
|c_1|^2\| \varphi \| ^2_2+\| \varphi \| _\infty \| q\circ \varphi
\| _2^2\leq |c_1|^2\| \varphi \| ^2_2+\| q\circ \varphi \|
_2^2\leq
\]
\[
(|c_1|^2_2+\| q\| _2^2)\| \varphi \| _2^2=\| \varphi \| _2^2\| p\|
^2_2.
\]
Because the polynomials in $H^2_0$ are dense in $H^2_0$, we get
$\|C_\varphi |H^2_0\| \le \|\varphi\|_2$ hence  $ \| C_\varphi
|H^2_0\| =\| \varphi \| _2$. \hfill \square
\bigskip

Condition (23) recalls \textit{Rudin's orthogonality condition}.
We say $\varphi $ satisfies \textit{Rudin's orthogonality
condition} if the family $\{ \varphi ^n:n=0,1,2,\dots \} $ is
orthogonal in $H^2$. It is easy to see that any $\varphi $ which
is the scalar multiple of an inner function fixing $0$ satisfies
the aforementioned condition. Recent results of Sundberg
~\cite{Sundberg}, show that there exist examples of symbols
satisfying \textit{Rudin's orthogonality condition}, other than
the multiples of inner functions fixing $0$, (thus answering a
question raised by Walter Rudin in 1988). Independently, Bishop
\cite{Bishop} obtained similar results using the pull--back
measure induced by $\varphi $. Obviously, the symbols fixing the
origin and satisfying (23) form a superset of those that satisfy
\textit{Rudin's orthogonality condition}. Elementary examples show
that the aforementioned superset is strictly larger: consider, for
instance, $\varphi (z)=(z^2+z^3)/2$.
\medskip

\textbf{Acknowledgement.} I am thankful to Paul Bourdon for
helping me prove Theorem 7 and his useful comments on the whole
paper.

\medskip

\noindent Department of Mathematics, University of Nebraska,
Omaha, NE 68182, USA. \textit{E-mail address:}
vmatache@mail.unomaha.edu

\end{document}